\documentclass[titlepage,11pt]{article}
\usepackage{latexsym}
\usepackage{amsfonts}
\usepackage{amsmath}
\newcommand\blackslug{\hbox{\hskip 1pt \vrule width 4pt height 8pt depth 1.5pt
        \hskip 1pt}}
\newcommand\bbox{\hfill \quad \blackslug \bigbreak}
\def\LL{,\ldots,}
\def\cupcup{\cup\cdots\cup}

\title{Polynomial bounds for chromatic number \\ VI. Adding a four-vertex path}
\author{Maria Chudnovsky\thanks{Supported by NSF DMS-EPSRC grant DMS-2120644.}\\
Princeton University, Princeton, NJ 08544
\\
\\
Alex Scott\thanks{Research supported by EPSRC grant EP/V007327/1.}\\
Mathematical Institute, University of Oxford, Oxford OX2 6GG, UK
\\
\\
Paul Seymour\thanks{Supported by AFOSR grants
A9550-19-1-0187 and FA9550-22-1-0234, and NSF grant DMS-2154169.}\\
Princeton University, Princeton, NJ 08544
\\
\\
Sophie Spirkl\thanks{We acknowledge the support of the Natural Sciences and Engineering Research
Council of Canada (NSERC), [funding reference number RGPIN-2020-03912].
Cette recherche a \'et\'e financ\'ee par le Conseil de recherches en sciences
naturelles et en g\'enie du Canada (CRSNG), [num\'ero de r\'ef\'erence
RGPIN-2020-03912].  }\\
University of Waterloo, Waterloo, Ontario N2L3G1, Canada}

\date{November 1, 2021; revised \today}

\newtheorem{thm}{}[section]

\newcommand{\Proof}{\noindent{\bf Proof.}\ \ }

\begin{document}
\maketitle
\begin{abstract}
A hereditary class of graphs is {\em $\chi$-bounded} if there is a function $f$ such that every graph $G$ in the class has chromatic number at most $f(\omega(G))$, where
$\omega(G)$ is the clique number of $G$; and the class is {\em polynomially $\chi$-bounded} if $f$ can be taken to be a polynomial.  The  Gy\'arf\'as-Sumner conjecture
asserts that, for every forest $H$, the class of $H$-free graphs (graphs with no induced copy of $H$) is $\chi$-bounded.    Let us say a forest $H$ is {\em good} if it satisfies the stronger property that the class of $H$-free graphs is polynomially $\chi$-bounded.

Very few forests are known to be good: for example, the goodness of the five-vertex path is open.
Indeed, it is not even known that if every component of a forest $H$ is good then $H$ is good, and in particular, it was not known that
        the disjoint union of two
        four-vertex paths is good. Here we show the latter (with corresponding polynomial $\omega(G)^{16}$); and more generally, that
        if $H$ is good then so is the disjoint union of $H$ and a four-vertex path.
                We also prove an even more general result: if every component of $H_1$ is good, and $H_2$ is any path (or broom)
then the class of graphs that are both $H_1$-free and $H_2$-free is polynomially $\chi$-bounded.
	
\end{abstract}

\section{Introduction}
A class of graphs is {\em hereditary} if it is closed under taking induced subgraphs; a hereditary class is
 {\em $\chi$-bounded} if there is a function $f$ such that every graph $G$ in the class has chromatic number at most $f(\omega(G))$, where
$\omega(G)$ is the clique number of $G$; and the class is {\em polynomially $\chi$-bounded} if $f$ can be taken to be a polynomial.
A graph is {\em $H$-free} if it has no induced subgraph isomorphic to $H$.

The Gy\'arf\'as-Sumner conjecture~\cite{gyarfas, sumner} asserts:
\begin{thm}\label{GSconj}
{\bf Conjecture: } For every forest $H$, the class of $H$-free graphs is $\chi$-bounded.
\end{thm}

There has been a great deal of recent progress on $\chi$-bounded classes (see \cite{survey} for a survey), although the Gy\'arf\'as-Sumner conjecture remains open.
In most cases, proofs of $\chi$-boundedness give fairly fast-growing functions, so it is interesting to ask:~when do we get the stronger property of polynomial $\chi$-boundedness?

A provocative conjecture of Louis Esperet~\cite{esperet} asserted that every $\chi$-bounded hereditary class is polynomially $\chi$-bounded,
but this was recently disproved by
Bria\'nski, Davies and Walczak \cite{BDW}.  So the question now is:~which hereditary classes are polynomially $\chi$-bounded?
In particular, can \ref{GSconj} be strengthened to polynomial $\chi$-boundedness? 
Let us say a graph $H$ is {\em good} if the class of $H$-free graphs is polynomially $\chi$-bounded. Perhaps every forest is good, but
the only trees currently known to be good are those not
containing the five-vertex path $P_5$~\cite{poly3}.
It is not known whether $P_5$ is good (although see~\cite{poly4} for the best 
current bounds for $H=P_5$; and see~\cite{poly5} for the case when $H$ is a general tree of radius two).

In the case of $\chi$-boundedness, it is not hard to show that a forest $H$ satisfies the Gy\'arf\'as-Sumner conjecture if and only if all its components do.  But
it has {\em not} been shown that if every component of a forest $H$ is good then $H$ is good.
Indeed, only some very restricted forests are known to be
good~\cite{Schiermeyer, poly2}.
One outstanding case was when $H$ is the forest $2P_4$, the disjoint union of two copies of
the four-vertex path $P_4$; and this was particularly annoying since the $P_4$-free graphs are very well-understood and
rather trivial.  We will prove that $2P_4$ is good, and indeed:
\begin{thm}\label{2P4}
If $G$ is $2P_4$-free, then $\chi(G)\le \omega(G)^{16}$.
\end{thm}

More generally, we will prove the following:
\begin{thm}\label{mainthm}
        If $H$ is a good forest, then the disjoint union of $H$
        and $P_4$ is also good.
\end{thm}

\ref{mainthm} is a consequence of the next result, about brooms.
A {\em $(k,d)$-broom} is a tree obtained from a $k$-vertex path with one end $v$ by adding $d$ new vertices adjacent to $v$, and
a {\em broom} is a tree that is a $(k,d)$-broom for some $k,d$. It is known that $(3,d)$-brooms are good~\cite{liu, poly3}, but
this is not known
for larger brooms (all of which contain $P_5$).  We will show the following, which implies \ref{mainthm}:
\begin{thm}\label{broomthm}
	Let $H_1$ be a forest such that every component of $H_1$ is good, and let $H_2$ be either
	a broom, or 
	the disjoint union of a good forest and a number of paths.
	Then there is a polynomial $\phi$ such that $\chi(G)\le \phi(\omega(G))$ for every $\{H_1,H_2\}$-free graph $G$.
\end{thm}
({\em $\{H_1,H_2\}$-free} means both $H_1$-free and $H_2$-free.)
To deduce \ref{mainthm} from \ref{broomthm}, let $H$ be a good forest, let $H_1=H_2$ be the disjoint union of $H$
	and $P_4$, and apply \ref{broomthm}.

	Some notation and terminology: if $G$ is a graph and $X\subseteq V(G)$, we denote by $G[X]$ the subgraph of $G$ induced on $X$, and
we sometimes write
$\chi(X)$ for $\chi(G[X])$ and $\omega(X)$ for $\omega(G[X])$.
Two disjoint subsets
$A,B\subseteq V(G)$ are {\em complete} if every vertex in $A$ is adjacent to every vertex of $B$, and 
{\em anticomplete} if there is no edge between $A,B$; and we say a vertex {\em $v$ is complete to $B$} if $\{v\}$ is complete to $B$, and so on.
A graph $G$ {\em contains} a graph $H$ if some
induced subgraph
of $G$ is isomorphic to $H$, and such a subgraph is a {\em copy} of $H$. The {\em cone} of a graph $H$
is obtained from $H$ by adding a new vertex adjacent to every vertex of $H$.

	Let us say a graph is {\em $0$-bad} if it is good; and a graph $J$ is $\beta$-bad, where $\beta\ge 1$ is an integer, if either
	$J$ is the disjoint union of two $(\beta-1)$-bad graphs, or $J$ is the cone of a $(\beta-1)$-bad graph, 
or $J$ is $(\beta-1)$-bad.
	In general, cones are not forests, so they are not good.
	Nevertheless, we will prove the following strengthening of \ref{broomthm}:
\begin{thm}\label{bigbrooms}
	Let $\beta\ge 0$, let $H_1$ be a $\beta$-bad graph, and let $H_2$ be either
        a broom, or
        the disjoint union of a good forest and a number of paths.
        Then there is a polynomial $\phi$ such that $\chi(G)\le \phi(\omega(G))$ for every $\{H_1,H_2\}$-free graph $G$.
\end{thm}

This implies several results that were previously known. For instance, in~\cite{Schiermeyer2} it is proved that:
\begin{thm}\label{pathandstuff}
        Let $H_1$ be either 
	\begin{itemize}
		\item the disjoint union of a complete graph and a good graph,  or
		\item the disjoint union of some complete graphs, or
		\item the cone of the disjoint union of some complete graphs.
	\end{itemize}
	Let $H_2$ be a path. Then
        there is a polynomial $\phi$ such that $\chi(G)\le \phi(\omega(G))$ for every $\{H_1,H_2\}$-free graph $G$.
\end{thm}
Some other results of~\cite{Schiermeyer2, Schiermeyer} are also special cases of \ref{bigbrooms}.

\section{Finding a disjoint union}

Suppose that $H$ is the disjoint union of good forests $H_1,H_2$.
Choose $c_1, c_2$ such that for $i=1,2$, every $H_i$-free graph $G$
satisfies $\chi(G)\le \omega(G)^{c_i}$.
Thus, if $G$ is $H$-free, 
we know that there do not exist disjoint, anticomplete subsets $P,Q\subseteq V(G)$ with 
$\chi(P)> \omega(P)^{c_1}$ and $\chi(Q)> \omega(Q)^{c_2}$; because then $G[P]$ is not $H_1$-free, and $G[Q]$ is not $H_2$-free,
and the union of a copy of $H_1$ in $G[P]$ and a copy of $H_2$ in $G[Q]$ gives a copy of $H$, which is impossible. 

But we do not 
really need $P,Q$ to be anticomplete.
It is enough that $\chi(P)>\omega(P)^{c_1}$, and $\chi(Q)>|H_1|r+ \omega(Q)^{c_2}$, where $r$ denotes the maximum over $v\in P$ 
of the chromatic number of the set of neighbours
of $v$ in $Q$; because then if we choose a copy $H_1'$ of $H$ in $G[P]$,
the chromatic number of the set of vertices in $Q$ with no neighbours in $V(H_1')$ is at least $\chi(Q)-|H_1|r>\omega(Q)^{c_2}$,
and so this set contains a copy of $H_2$, a contradiction. In the proof to come later in the paper, this is the only way we will ever 
use that $G$ is $H$-free; and so we might as well
prove a stronger theorem, replacing the hypothesis that $G$ is $H$-free with the weaker hypothesis that there is no
suitable pair $(P,Q)$ in $G$.
 
Thus we will be excluding pairs of disjoint sets $P,Q$ where $\chi(P)$
is at least some power of $\omega(P)$, and for each vertex in $P$, its set of neighbours in $Q$ has chromatic number at most some
$r$ that is small
compared with the chromatic number of $Q$. 

In our proof, it happens that when we find a suitable pair $(P,Q)$, it comes equipped with an extra vertex $v$ that is 
complete to $P$ and anticomplete to $Q$; so we might as well prove that there is a ``suitable triple'' $(v,P,Q)$. Such a thing 
will also allow us to handle cones.

We denote the set of nonnegative integers by $\mathbb{N}$, and say a function $\phi:\mathbb{N}\rightarrow\mathbb{N}$ is 
{\em non-decreasing} if
$\phi(x)\le \phi(x')$ for all $x,x'\in \mathbb{N}$ with $x\le x'$. 

Let $\psi:\mathbb{N}\rightarrow\mathbb{N}$ be non-decreasing, and let $q\ge 0$ be an integer.
We say a {\em $(\psi,q)$-scattering} in a graph $G$ is a triple $(v,P,Q)$ where:
\begin{itemize}
	\item $P,Q$ are disjoint subsets of $V(G)$, and $v\in V(G)\setminus (P\cup Q)$;
	\item $\{v\}$ is complete to $P$ and anticomplete to $Q$;
	\item $\chi(P)> \psi(\omega(P))$; and
	\item $\chi(Q)> qr+\psi(\omega(Q))$, where $r$ is the maximum, over $u\in P$, of the chromatic number of 
		the set of neighbours of $u$ in $Q$.
\end{itemize}

Thus we will replace the hypothesis in \ref{bigbrooms} that $G$ is $H_1$-free and $H_1$ is $\beta$-bad,
with the hypothesis that $G$ contains no $(\psi,q)$-scattering, for appropriate $\psi,q$. We will show:
\begin{thm}\label{scattering}
	Let $\psi:\mathbb{N}\rightarrow\mathbb{N}$ 
	be a non-decreasing polynomial and let $q\in \mathbb{N}$. 
        Let $H_2$ be either a broom, or the disjoint union of a good forest and a number of paths.
	Then there is a polynomial $\phi:\mathbb{N}\rightarrow\mathbb{N}$ such that if $\chi(G)> \phi(\omega(G))$ and $G$ contains no 
	$(\psi,q)$-scattering, then $G$ contains $H_2$.
\end{thm}
\noindent{\bf Proof of \ref{bigbrooms}, assuming \ref{scattering}.\ }
We proceed by induction on $\beta$.
Let $H_1$ be $\beta$-bad, and let $H_2$ be either 
a broom, or
the disjoint union of a good forest and a number of paths.

If $H_1$ is good, the result is true, so we assume that $H_1$ is not good, and therefore $\beta\ge 1$.
Thus either $H_1$ is the disjoint union of two $(\beta-1)$-bad graphs $J_1,J_2$, or the cone of a $(\beta-1)$-bad graph $J_1$
(and in this case let $J_2$ be the null graph).
From the inductive hypothesis on $\beta$, for $i = 1,2$ there is a non-decreasing polynomial $\phi_i$ such that if $G$ is $H_2$-free 
and $J_i$-free
then $\chi(G)\le \phi_i(\omega(G))$, and by replacing $\phi_1,\phi_2$ by $\phi_1+\phi_2$ we may assume that $\phi_1=\phi_2$. 

Let $q=|J_1|$. By \ref{scattering}, there is a non-decreasing polynomial $\phi$ such that 
if $\chi(G)> \phi(\omega(G))$ and contains no $(\phi_1,q)$-scattering, then $G$ contains $H_2$.
We claim that $\phi$ satisfies \ref{bigbrooms}.

Let $G$ be $\{H_1, H_2\}$-free, and suppose that $\chi(G)> \phi(\omega(G))$. Since $G$ is $H_2$-free, it follows from the choice of $\phi$
that $G$ contains a 
$(\phi_1,q)$-scattering $(w,P,Q)$ say.
Let $r$ be the maximum, over $v\in P$, of the chromatic number of
the set of neighbours of $v$ in $Q$. Since $\chi(P)> \phi_1(\omega(P))$, there is an induced subgraph of $G[P]$ isomorphic to $J_1$,
say $J_1'$. Hence $G$ contains the cone of $J_1$, so we may assume that $H_1$ is the disjoint union of $J_1,J_2$. 
The set of vertices in $Q$ with a neighbour in $V(J_1')$ has chromatic number at most $r|J_1|$, and since
$$\chi(Q)> |J_1|r+\phi_2(\omega(Q)),$$
it follows that the set (say $Q'$) of vertices in $Q$ that are anticomplete to $J_1'$ has chromatic number more than 
$\phi_2(\omega(Q))$. From the choice of $\phi_2$, and since $G$ is $H_2$-free, it follows that $G[Q']$ is not $J_2$-free; but then, 
combining
this copy of $J_2$ with $J_1'$, we find a copy of $H_1$ in $G$, a contradiction. This proves \ref{bigbrooms}.~\bbox

Let $\sigma:\mathbb{N}\rightarrow\mathbb{N}$  be a non-decreasing function. We say a subgraph $P$ of a graph $G$
is {\em $\sigma$-nondominating} if there is a set $X\subseteq V(G)\setminus V(P)$, anticomplete to $V(P)$,  with $\chi(X)> \sigma(\omega(X))$.
Next we will show that to prove \ref{scattering} it suffices to prove the following:

\begin{thm}\label{sweep}
        Let $\psi,\sigma:\mathbb{N}\rightarrow\mathbb{N}$ be non-decreasing polynomials, and let $q\ge 0$ an integer.
        Let $H$ be a broom, and let $J$ be a path.
        Then there is a non-decreasing polynomial $\phi:\mathbb{N}\rightarrow\mathbb{N}$ such that if $G$ is a graph, and 
	$\chi(G)> \phi(\omega(G))$, and 
	$G$ contains no $(\psi,q)$-scattering, then $G$ contains $H$ and a $\sigma$-nondominating copy of $J$.
\end{thm}
\noindent{\bf Proof of \ref{scattering}, assuming \ref{sweep}.\ }
Let $\psi,q, H_2$ be as in \ref{scattering}. If $H_2$ is a broom, then \ref{scattering} follows immediately from \ref{sweep} 
(setting $H=H_2$ and setting $J$ to be some path, for instance the one-vertex path). Thus we assume that 
$H_2$ is the disjoint union of a good forest 
$J_1$ and a forest $J_2$ that is a disjoint union of paths. Let
$\sigma:\mathbb{N}\rightarrow\mathbb{N}$ be a non-decreasing function such that every $J_1$-free graph $G$ has chromatic number at most 
$\sigma(\omega(G))$;
and choose a path $J$ such that $J_2$ is an induced subgraph of $J$.
By \ref{sweep} (setting $H$ to be some broom, for instance with one vertex) there is a non-decreasing polynomial 
$\phi:\mathbb{N}\rightarrow\mathbb{N}$ such that if $\chi(G)> \phi(\omega(G))$ 
and $G$ contains no $(\psi,q)$-scattering, then
$G$ contains a $\sigma$-nondominating copy $J'$ of $J$.

We claim that $\phi$ satisfies \ref{scattering}.
Thus we must show that if $G$ is $H_2$-free and contains no $(\psi,q)$-scattering then $\chi(G)\le \phi(\omega(G))$.
Suppose not. By the choice of $f$, and since $G$ contains no $(\psi,q)$-scattering, it follows that 
$G$ contains a copy $J'$ of $J$, such that there is a set $X\subseteq V(G)$ with $\chi(X)> \sigma(\omega(X))$
anticomplete to $V(J_1')$. But since $\chi(X)> \sigma(\omega(X))$, it follows that
$G[X]$ contains $J_1$, and since $J$ contains $J_2$, and $V(J)$
is anticomplete to $X$, it follows that $G$ contains $H_2$. This proves \ref{scattering}.~\bbox

We remark that there is an appealing possible strengthening of \ref{sweep}, that we could not prove:
\begin{thm}\label{bettersweep}
        {\bf Conjecture: }Let $\psi,\sigma:\mathbb{N}\rightarrow\mathbb{N}$ be non-decreasing polynomials, let $q\ge 0$ an integer, and let $H$ be a broom.
        Then there is a non-decreasing polynomial $\phi:\mathbb{N}\rightarrow\mathbb{N}$ such that if $G$ is a graph, and
        $\chi(G)> \phi(\omega(G))$, and
        $G$ contains no $(\psi,q)$-scattering, then $G$ a $\sigma$-nondominating copy of $H$.
\end{thm}

Let us say a graph $H$ is {\em self-isolating} if for every non-decreasing
polynomial $\psi:\mathbb{N}\rightarrow \mathbb{N}$, there is a polynomial
$\phi:\mathbb{N}\rightarrow \mathbb{N}$ with the following property: for every graph $G$ with $\chi(G)> \phi(\omega(G))$,
there exists $A\subseteq V(G)$ with $\chi(A)>\psi(\omega(A))$, such that either
\begin{itemize}
        \item $G[A]$ is $H$-free, or
        \item $G$ contains a copy $H'$ of $H$ such that $V(H')$ is disjoint from and anticomplete to $A$.
\end{itemize}
Which graphs are self-isolating? It is proved in~\cite{poly2} that stars are self-isolating, and we will show in~\cite{poly7}
that complete graphs and complete bipartite graphs are self-isolating. Let us observe that \ref{sweep} implies that:
\begin{thm}\label{pathisolation}
Every path is self-isolating.
\end{thm}
\Proof
Let $J$ be a path, and let $\psi:\mathbb{N}\rightarrow \mathbb{N}$ be a non-decreasing polynomial. 
Choose $\phi$ satisfying \ref{sweep} with $H=J$ and $\sigma=\psi$ and $q=|J|$,
and let $G$ be a graph with $\chi(G)> \phi(\omega(G))$.
We claim that either there is a $\psi$-nondominating copy of $J$ in $G$, or
there exists $A\subseteq V(G)$ with $\chi(A)>\psi(\omega(A))$ such that $G[A]$ is $J$-free.
By \ref{sweep} we may assume that there is a  $(\psi, q)$-scattering $(w, P,Q)$ in $G$. If $G[P]$ is $J$-free, the claim holds,
so we assume that there is a copy $J'$ of $J$ in $G[P]$. Thus $|J'|=q$.
Let $r$ be the maximum over $v\in P$ of the chromatic number of the
set of neighbours of $v$ in $Q$.
The set of vertices in $Q$ with a neighbour in $V(J')$ has chromatic number at most $|J'|r=qr$; and $\chi(Q)> \psi(\omega(Q))+qr$
from the definition of a $(\psi, q)$-scattering. Consequently $J'$ is $\psi$-nondominating, and hence $J$ is self-isolating. 
This proves \ref{pathisolation}. ~\bbox

\section{Constructing a horn}
Let $d\ge 0$ be an integer. If $A,B\subseteq V(G)$ are disjoint, we say that $A$ is {\em $d$-dense} to $B$
if
for every vertex $v\in A$, the set of non-neighbours of $v$ in $B$ has chromatic number at most $d$.
Let us say a {\em $(d,z)$-horn} in a graph $G$ is a triple $(v,A,B)$ where
\begin{itemize}
        \item $A,B$ are disjoint subsets of $V(G)$, and $v\in V(G)\setminus (A\cup B)$;
        \item $v$ is complete to $A$ and anticomplete to $B$; and 
        \item there is no $Z\subseteq A\cup B$ with $\chi(Z)\le z$ such that $A\setminus Z$ is $d$-dense to $B\setminus Z$.
\end{itemize}

We will need a $(d,z)$-horn $(v,A,B)$ where $z$ is  at least some large function of the clique number of $A\cup B$, and 
this section produces such a horn.
We show in \ref{getbigcone} that if $G$ has sufficiently large chromatic number (and, for convenience, all its proper induced subgraphs have
smaller chromatic number), then either $G$ contains both a $(k,s)$-broom and a $\sigma$-nondominating $k$-vertex path, or
$G$ contains a $(d,z)$-horn. To complete the proof of \ref{sweep}, it therefore suffices to handle graphs $G$ that contain
$(d,z)$-horns, for suitably chosen values of $d,z$, and we will do so in the next section.

We will use the following well-known version of Ramsey's theorem, proved (for instance)
in~\cite{poly2} ($|G|$ denotes the number of vertices of $G$):
\begin{thm}\label{ramsey}
        Let $x\ge 2$ and $y\ge 1$ be integers. For a graph $G$, if $|G|\ge x^y$, then $G$ has either a clique of
        cardinality $x+1$, or a stable set of cardinality $y$.
\end{thm}

If $v\in V(G)$, we denote by $N(v)$ or $N_G(v)$ the set of all neighbours of $v$ in $G$. First, we need a result of 
Gy\'arf\'as~\cite{gyarfasworld} (we give the well-known proof, because it is so pretty.)
\begin{thm}\label{pathnbrs}
	Let $k\ge 1$ and $x\ge 0$ be integers.
	Let $G$ be a connected graph such that $\chi(N(v))\le  x$ for every vertex $v$.
	Let $H$ be a connected induced subgraph of $G$, and let $v\in V(G)\setminus V(H)$
	with a neighbour in $V(H)$. If $\chi(H)>  (k-2)x$, there is an induced $k$-vertex path of $G$ with one end $v$
	and all other vertices in $V(H)$.
\end{thm}
\Proof We proceed by induction on $k$. The result is clear for $k\le 2$, so we assume that $k\ge 3$.
Let $J$ be obtained from $H$ by deleting all vertices in $N(v)$; thus $\chi(J)>(k-3)x>0$, and so
there is a component $H'$ of $J$ with chromatic number more than $(k-3)x$. Let $v'\in N(v)\cap V(H)$
with a neighbour in $V(H')$. From the inductive hypothesis applied to $v',H'$, there is an induced $(k-1)$-vertex path 
of $G$ with one end $v'$
and all other vertices in $V(H')$. Appending $v$ to this path proves \ref{pathnbrs}.~\bbox

We deduce:
\begin{thm}\label{bigpath}
	Let $\sigma:\mathbb{N}\rightarrow\mathbb{N}$ be non-decreasing, let $k,x\ge 1$
	be integers, and let $G$ be a graph. If $\chi(N(v))\le  x$ for every $v\in V(G)$, and $\chi(G)>kx+\sigma(\omega(G))$,
	then there is a $\sigma$-nondominating $k$-vertex induced path $P$ in $G$.
\end{thm}
\Proof We may assume that $G$ is connected; choose $v\in V(G)$. Since $\chi(G\setminus v)> kx-1\ge (k-2)x$, \ref{pathnbrs}
(applied to $v$ and to a component of $G\setminus v$ of maximum chromatic number) implies that $G$ contains a 
$k$-vertex induced path $P$.
The set of vertices of $G$ with a neighbour in $V(P)$ has chromatic number at most $kx$, and the result follows. This proves
\ref{bigpath}.~\bbox

The next result is also essentially due to Gy\'arf\'as (mentioned  in~\cite{gyarfasworld}):
\begin{thm}\label{broomnbrs}
	Let $H$ be a $(k,s)$-broom, and suppose that $G$ is $H$-free, and $\chi(N(v))\le  x$ for every $v\in V(G)$. Then 
	$$\chi(G)\le \max(\omega(G)^{2s}, (2s+1)(x+1)+(k-2)x).$$
\end{thm}
\Proof
Suppose that $\chi(G)>\max(\omega(G)^{2s}, (2s+1)(t+1)+(k-2)x)$. We may assume that $G$ is connected. 
If every vertex of $G$ has degree less than $\omega(G)^{2s}$
then $\chi(G)\le \omega(G)^{2s}$, a contradiction, so some vertex $v$ has at least $\omega(G)^{2s}$ neighbours. By
\ref{ramsey} applied to $G[N(v)]$, there is a stable set $S$ of neighbours of $v$, with $|S|=2s$.
Let $M$ be the set of all vertices of $G$ that do not belong to $S\cup \{v\}$ and have a neighbour in $S\cup \{v\}$.
Thus $\chi(M)\le (2s+1)x$. Let $H$ be a component of $G\setminus (M\cup S\cup \{v\})$ of maximum chromatic number;
then $\chi(H)\ge \chi(G)-(2s+1)(x+1)> (k-2)x$. Choose $u\in M\cup S\cup \{v\}$ with a neighbour in $V(H)$. Since no vertex 
of $S\cup \{v\}$ has a neighbour in $V(H)$, from the definition of $M$, it follows that $u\in M$.
By \ref{pathnbrs} applied to
$u,H$, there is an induced $k$-vertex path $P$ of $G$ with one end $u$
and all other vertices in $V(H)$. Thus $u$ is the only vertex of $P$ with a neighbour in $S\cup \{v\}$. 
If $u$ is adjacent to at least $s$ vertices in $S$, then the subgraph induced on $V(P)$ and some $s$ of these neighbours is a 
$(k,s)$-broom, a contradiction. Thus there exists $S'\subseteq S$ with $|S'|=s$, such that all vertices in $S'$ are nonadjacent to $u$.
If $u$ is adjacent to $v$, the subgraph induced on $V(P)\cup S\cup \{v\}$ is a
$(k+1,s)$-broom, a contradiction. Thus $u$ is adjacent to some $w\in S\setminus S'$, and nonadjacent to $v$. But then the 
subgraph induced on $V(P)\cup S'\cup \{v,w\}$ is a $(k+2,s)$-broom, a contradiction. This proves \ref{broomnbrs}.~\bbox

\begin{thm}\label{getbigcone}
	Let $\sigma:\mathbb{N}\rightarrow\mathbb{N}$ be non-decreasing. Let $k,s,d,z\ge 0$ and $c\ge 2s$ be integers.
	Let $G$ be a graph such that
	\begin{align*}
		\chi(G)&>\omega(G)^c;\\
		 \chi(G')&\le \omega(G')^c \text{ for every induced subgraph $G'$  of   $G$ with $G'\ne G$;}\\
		\omega(G)^c &\ge (\omega(G)-1)^c +z+d\omega(G)+2;\\
		\omega(G)^c &\ge  (2s+1)(z+1)+(k-2)z; \text{ and }\\
		\omega(G)^c &\ge kz+\sigma(\omega(G)).
	\end{align*}
Then either 
	\begin{itemize}
		\item $G$ contains a $(d,z)$-horn; or
		\item $G$ contains a $(k,s)$-broom, and a $\sigma$-nondominating $k$-vertex path.
        \end{itemize}
\end{thm}
\Proof
Suppose that $\chi(N(v))\le z$ for every vertex $v\in V(G)$. By \ref{broomnbrs}, and since
$$\chi(G)>\omega(G)^c\ge \max(\omega(G)^{2s}, (2s+1)(z+1)+(k-2)z)$$
(because $c\ge 2s$), it follows that $G$ contains a $(k,s)$-broom.
By \ref{bigpath}, since 
$\chi(G)-kz>\sigma(\omega(G))$,
there is a $\sigma$-nondominating $k$-vertex induced path $P$ in $G$, 
	and so the second bullet holds.

	Thus we assume that $\chi(N(v))> z$ for some vertex $v$. Let $A$ be the set of neighbours of $v$, and 
	$B=V(G)\setminus (A\cup \{v\})$. We claim that $(v,A,B)$ is a $(d,z)$-horn. Suppose not; then there exists
	$Z\subseteq A\cup B$ with $\chi(Z)\le z$, such that $A\setminus Z$ is $d$-dense to $B\setminus Z$. Let 
	$P\subseteq A\setminus Z$ be a clique with cardinality $p=\omega(A\setminus Z)$.
	Then $p\ge 1$, since $\chi(Z)\le z<\chi(A)$; and $p<\omega(G)$ since otherwise adding $v$ would give a clique of cardinality
	$\omega(G)+1$. For each $u\in P$, the set of vertices in $B\setminus Z$ nonadjacent to $u$ has chromatic number at most $d$,
	since $A\setminus Z$ is $d$-dense to $B\setminus Z$; and so the set of vertices in $B$ with a non-neighbour in $P$
	has chromatic number at most $pd\le d\omega(G)$. The set of vertices in $B$ complete to $P$ has clique number at most $\omega(G)-p$
	and so has chromatic number at most $(\omega-p)^c$. Hence $\chi(B\setminus Z)\le pd+(\omega(G)-p)^c$, and so
	$$\chi(G)\le \chi(Z)+\chi(A\setminus Z)+\chi(B\setminus Z)+1\le z+p^c+ d\omega(G)+(\omega(G)-p)^c+1.$$
	Since $1\le p\le \omega(G)-1$, $p^c+(\omega(G)-p)^c\le (\omega(G)-1)^c+1$, and so
	$$\omega(G)^c<\chi(G)\le  z+ d\omega(G)+(\omega(G)-1)^c+2,$$
	a contradiction. This proves \ref{getbigcone}.~\bbox


\section{Making taller horns}

In this section we prove \ref{sweep}, and hence complete the proofs of \ref{scattering}, \ref{bigbrooms}, \ref{broomthm}, and therefore
\ref{mainthm}. Because of \ref{getbigcone}, we may assume that $G$ contains a $(d,z)$-horn, for some suitable values of $d,z$;
and now we will show that, provided that $G$ does not contain the proscribed scattering,  we can use this horn 
to make a ``$k$-tall'' $(d',z')$-horn, which is a horn with a $k$-vertex path
appended to its distinguished vertex. From such a horn, it is easy to obtain a
$(k,s)$-broom and a $\sigma$-nondominating $k$-vertex path, to satisfy \ref{sweep}. The main step is therefore to convert an $\ell$-tall
horn to an $(\ell+1)$-tall horn, and for that we need the next result.

If $d,z,\omega\ge 0$ are integers, a graph $G$ is {\em $(d,z, \omega)$-unsplittable} if there is no partition $(A,B,Z)$ of $V(G)$
such that $\chi(Z)\le z$, and $\chi(A), \chi(B)> d\omega$, and $A$ is $d$-dense to $B$.
We begin with:
\begin{thm}\label{splitpartition}
	If $d,z\ge 0$ are integers, every graph $G$ admits a partition $(D_0,D_1\LL D_k)$  of its 
	vertex set with $k\le \omega(G)$ 
	such that $\chi(D_0)\le z\omega(G)$ and $G[D_i]$
	is $(d,z,\omega(G))$-unsplittable for $1\le i\le k$.
\end{thm}
\Proof
We may assume that $G$ is not $(d,z, \omega(G))$-unsplittable, and so it admits a partition $(D_0,D_1,D_2)$ such that $\chi(D_0)\le z$,
$\chi(D_1), \chi(D_2)> d\omega(G)$, and $D_1$ is $d$-dense to $D_2$. Hence we may choose
$k\ge 2$ maximum such that there is a sequence $D_0,D_1\LL D_k$ of pairwise disjoint subsets of $V(G)$ with union $V(G)$,
and with the following properties:
\begin{itemize}
	\item $\chi(D_0)\le (k-1)z$
        \item $D_i$ is $d$-dense to $D_j$ for $1\le i<j\le k$; and 
	\item $\chi(D_i)> d\omega(G)$ for $1\le i\le k$.
\end{itemize}
We claim:
\\
\\
(1) {\em $k\le \omega(G)$.}
\\
\\
Suppose that $k>\omega(G)$, and define $d_i\in D_i$ for $1\le i\le \omega(G)+1$ inductively as follows. 
Let $1\le i\le \omega(G)+1$, 
and suppose that $d_1\LL d_{i-1}$
have been defined, all pairwise adjacent. 
The set of vertices in $D_i$ that have a non-neighbour among $d_1\LL d_{i-1}$ has chromatic number at most
$$(i-1)d\le d\omega(G) <\chi(D_i),$$
and so some vertex $d_i\in D_i$ is adjacent to all of $d_1\LL d_{i-1}$. This completes the inductive
definition. But then $\{d_1\LL d_{\omega(G)+1}\}$ is a clique of $G$, contradicting the definition of $\omega(G)$. This proves (1).
\\
\\
(2) {\em For $1\le i\le k$, $G[D_i]$ is $(d,z,\omega(G))$-unsplittable.}
\\
\\
Suppose that $(A,B, Z)$ is a partition of $D_i$
such that $\chi(Z)\le z$, and $\chi(A), \chi(B)> d\omega(G)$, and $A$ is $d$-dense to $B$.
Then the sequence
$$(D_0\cup Z, D_1\LL D_{i-1}, A, B, D_{i+1}\LL D_k)$$
contradicts the maximality of $k$. This proves (2).

\bigskip

From (1), (2), this proves \ref{splitpartition}.~\bbox

Let $(v,A,B)$ be a $(d,z)$-horn in a graph $G$, and let $k\ge 1$ be an integer. We say that $(v,A,B)$ is {\em $k$-tall} if
there is an induced path $R$ in $G$ with $k$ vertices, with one end $v$, such that $V(R)\setminus \{v\}$ is disjoint from and anticomplete
to $A\cup B$. Thus every  $(d,z)$-horn is 1-tall.
We use \ref{splitpartition} to prove a result which is the heart of the paper:
\begin{thm}\label{findmixed}
	Let $G$ be a graph, let $d,z,d',z',q\ge 0$ be integers, and let $\psi:\mathbb{N}\rightarrow\mathbb{N}$
	be non-decreasing, satisfying:
	\begin{align*}
		z&\ge \left(2\psi(\omega(G))+(1+q)z'+qd'\omega(G) \right)\omega(G)\\
		d&\ge \left(z'+d'\omega(G)\right)\omega(G).
	\end{align*}
	Let $(v,A,B)$ be an $\ell$-tall $(d,z)$-horn in a graph $G$, for some $\ell\ge 1$. Then either
	\begin{itemize}
		\item there exist $P\subseteq A$ and $Q\subseteq B$ such that $(v,P,Q)$ is a $(\psi,q)$-scattering; or
		\item there exist $v'\in A$ and disjoint subsets $A',B'$ of $B$ such that $(v', A', B')$ is an $(\ell+1)$-tall
			$(d',z')$-horn.
	\end{itemize}
\end{thm}
\Proof Let $p=\psi(\omega(G))$. By \ref{splitpartition}, $B$ admits a partition $(D_0,D_1\LL D_k)$ with 
$k\le \omega(G)$ such that $\chi(D_0)\le z'\omega(G)$ and $G[D_i]$ is $(d',z',\omega(G))$-unsplittable for $1\le i\le k$.
For $1\le i\le k$, if $\chi(D_i)\le  q(z'+d'\omega(G))+p$ let $P_i=\emptyset$, and if 
$\chi(D_i)> q(z'+d'\omega(G))+p$ let 
$P_i$ be the set of vertices $a\in A$ 
such that $\chi(U)\le  z'+d'\omega(G)$, where $U$ is the set 
of neighbours of $a$ in $D_i$. Let $P= P_1\cupcup P_k$.

Suppose that $\chi(P_i)> p$, for some $i\in \{1\LL k\}$. Consequently $P_i\ne \emptyset$, and so 
$$\chi(D_i)> q(z'+d'\omega(G))+p\ge q(z'+d'\omega(G))+ \psi(\omega(D_i));$$
and for each $a\in P_i$, $\chi(U)\le  z'+d'\omega(G)$, where $U$ is the set
of neighbours of $a$ in $D_i$. It follows that $(v,P_i,D_i)$ is a $(\psi,q)$-scattering
and the first bullet of the theorem holds. 
Thus we may assume that $\chi(P_i)\le  p$ for $1\le i\le k$, 
and consequently $\chi(P)\le p\omega(G)$.

Let $Z$ be the union of $P, D_0$, and all the sets $D_i$ with $1\le i\le k$ such that 
$$\chi(D_i)\le q(z'+d'\omega(G))+p.$$
Consequently 
$$\chi(Z)\le 2p\omega(G)+ z'\omega(G) + q(z'+d'\omega(G))\omega(G)\le z.$$
Since $(v,A,B)$ is a $(d,z)$-horn, it follows that $A\setminus Z$ is not $d$-dense to $B\setminus Z$; and so there exists 
$v'\in A\setminus P$ such that the set of vertices in $B\setminus Z$ that are nonadjacent to $v'$ has chromatic number
more than $d$. Since $B\setminus Z$ is the union of the sets $D_i$ with $\chi(D_i)> q(z'+d'\omega(G))+p$, 
there exists $i\in \{1\LL k\}$
with $\chi(D_i)\ge q(z'+d'\omega(G))+p$ such that the set $B'$ of vertices in $D_i$ nonadjacent to $v'$ has chromatic number more 
than $d/\omega(G)$. Since
$v'\notin P$, the set $A'$ of neighbours of $v'$ in $D_i$ has chromatic number more than $d'\omega(G)+z'$. 

Let $Z'\subseteq D_i$ with $\chi(Z')\le z'$. Thus $\chi(A'\setminus Z')\ge \chi(A')-\chi(Z')> d'\omega(G)$; and
$\chi(B'\setminus Z')> d/\omega(G)-z'\ge d'\omega(G)$. Since $G[D_i]$ is $(d',z',\omega(G))$-unsplittable, it follows
that $A'\setminus Z'$ is not $d'$-dense to $B'\setminus Z'$. This proves that  $(v', A', B')$ is 
a $(d', z')$-horn. 

Since $(v,A,B)$ is $\ell$-tall, there is an $\ell$-vertex induced path $R$ of $G$ with one end $v$, such that $V(R)\setminus \{v\}$
is disjoint from and anticomplete to $A\cup B$. Then $R'=G[V(R)\cup \{v'\}]$ is an $(\ell+1)$-vertex path, and since $V(R)$
is anticomplete to $B$ and hence to $A'\cup B'$, it follows that $(v', A', B')$ is $(\ell+1)$-tall, and so the second bullet of the theorem holds. This proves \ref{findmixed}.~\bbox

Now we prove \ref{sweep}, which we restate:
\begin{thm}\label{sweep2}
        Let $k,s\ge 1$ and $q\ge 0$ be integers, and let $\psi,\sigma:\mathbb{N}\rightarrow\mathbb{N}$
be non-decreasing polynomials.
        Then there exists an integer $c\ge 0$ such that if $G$ is a graph with 
        $\chi(G)> \omega(G)^c$, and
	$G$ contains no $(\psi,q)$-scattering, then $G$ contains a $(k,s)$-broom and a $\sigma$-nondominating $k$-vertex path.
\end{thm}
\Proof
Let $\zeta_k:\mathbb{N}\rightarrow\mathbb{N}$ be the polynomial defined by $\zeta_k(x)=\sigma(x)+x^s$, and let $\delta_k(x)=0$. 
For $i=k-1\LL 1$, define polynomials $\zeta_i, \delta_i:\mathbb{N}\rightarrow\mathbb{N}$ by
\begin{align*}
	\zeta_i(x)&=2x\psi(x)+(1+q)x\zeta_{i+1}(x)+qx^2\delta_{i+1}(x)\\
	\delta_i(x)&= x\zeta_{i+1}(x)+x^2\delta_{i+1}(x).
\end{align*}
Choose an integer $c\ge 2s$ such that
\begin{align*}
	x^c &\ge (x-1)^c +\zeta_1(x)+x\delta_1(x)+2\\
		x^c &\ge  (2s+1)(\zeta_1(x)+1)+(k-2)\zeta_1(x), \text{ and }\\
		x^c &\ge k\zeta_1(x)+\sigma(x)
\end{align*}
for all integers $x\ge 2$. We claim that $c$ satisfies \ref{sweep2}. To see this, let 
$G$ be a graph with
        $\chi(G)> \omega(G)^c$, and suppose that 
        $G$ contains no $(\psi,q)$-scattering. We must show that $G$ contains a $(k,s)$-broom and a $\sigma$-nondominating $k$-vertex path.
We show this by induction on $|G|$. If there is 
an induced subgraph $G'$ of $G$ with $G'\ne G$ and $\chi(G')>\omega(G')^c$, then $G'$ contains no $(\psi,q)$-scattering, and from 
the inductive hypothesis, $G'$ contains a $(k,s)$-broom and a $\sigma$-nondominating $k$-vertex path, and hence so does $G$, as required.
We may assume then that there is no such $G'$. Since $\chi(G)> \omega(G)^c$, it follows that $\omega(G)\ge 2$, and so
the five displayed inequalities of \ref{getbigcone} hold with $z,d$ replaced by $\zeta_1(\omega(G)), \delta_1(\omega(G))$ respectively.
From \ref{getbigcone}, we may assume that $G$ contains a $(\delta_1(\omega(G)),\zeta_1(\omega(G)))$-horn, which is therefore 1-tall.

From \ref{findmixed},
it follows that for $i=2\LL k$, $G$ contains an $i$-tall $(\delta_i(\omega(G)), \zeta_i(\omega(G)))$-horn, and so contains a $k$-tall
$(0, z)$-horn $(v,A,B)$ say, where $z=\zeta_k(\omega(G))$. 
Since this horn is $k$-tall, there is a
$k$-vertex induced path $R$ of $G$ with one end $v$, such that
$V(R)\setminus \{v\}$ is disjoint from and anticomplete to $A\cup B$. From the definition of a $(0,z)$-horn, $\chi(A),\chi(B)> z$.
Since $\chi(A)> z\ge \omega(A)^s$,
\ref{ramsey} implies that there is a stable set $S\subseteq A$ with $|S|=s$, and so $G[V(R)\cup S]$ is a $(k,s)$-broom.
Since $\chi(B)> z> \sigma(\omega(B))$, and $V(R)$ is anticomplete to $B$, $R$ is $\sigma$-nondominating.
This proves \ref{sweep2}.~\bbox

Finally, we will go through the calculations of the proof of \ref{sweep2}, to prove \ref{2P4}, which we restate:
\begin{thm}\label{2P42}
If $G$ is $2P_4$-free, then $\chi(G)\le \omega(G)^{16}$.
\end{thm}
\Proof
Let $\alpha$ be the polynomial where $\alpha(x)=x$ for all $x$. If $G$ is $2P_4$-free, then $G$ contains no $(\alpha,4)$-scattering,
and contains no $\alpha$-nondominating $4$-vertex path; so we will follow the proof of \ref{sweep2}, taking
$k=q=4$, $s=1$, and $\psi=\sigma=\alpha$.

Thus, from the definitions, we have that for all $x\ge 0$:
\begin{align*}
        \zeta_4(x) &=2x\\
        \delta_4(x) &=0\\
        \zeta_3(x)&=12x^2\\
        \delta_3(x)&= 2x^2\\
        \zeta_2(x)&=2x^2+60x^3+8x^4\\
        \delta_2(x)&= 12x^3+2x^4\\
        \zeta_1(x)&=2x^2+ 10x^3+300x^4 + 88x^5+ 8x^6\\
        \delta_1(x)&= 2x^3+60x^4+ 20x^5+ 2x^6.
\end{align*}
Then we must choose an integer $c\ge 2$ such that
\begin{align*}
        x^c &\ge (x-1)^c + 2+ 2x^2+ 10x^3+302x^4 + 148x^5+ 28x^6 +2x^7\\
                x^c &\ge  3+10x^2+ 50x^3+ 1500x^4+ 440x^5+ 40x^6, \text{ and }\\
                x^c &\ge x+8x^2+ 40x^3+1200x^4 + 352x^5+ 32x^6.
\end{align*}
for all integers $x\ge 2$. Thus we may take $c=16$. This proves \ref{2P42}.~\bbox


\end{document}